\documentclass[11pt,a4paper]{amsart}
\usepackage{amssymb,amsmath}
\numberwithin{equation}{section}
\vspace{8cm}
\newtheorem{theorem}{Theorem}[section]
\newtheorem{proposition}[theorem]{Proposition}
\newtheorem{corollary}[theorem]{Corollary}
\newtheorem{lemma}[theorem]{Lemma}

\newtheorem{remark}[theorem]{Remark}

\newcommand{\R}{{\mathbb R}}
\newcommand{\Rn}{{\mathbb R}^n}

\def\d0{\delta_0}

\def\sc0{(x_0,y_0)}
\newcommand{\beq}{\begin{equation}}
\newcommand{\eeq}{\end{equation}}
\newcommand{\refe}[1]{{\rm (\ref{#1})}}
\newcommand{\dis}{\displaystyle}

\newcommand{\Mm}{\mathcal{M}^-_{\lambda, \Lambda}}
\newcommand{\Mp}{\mathcal{M}^+_{\lambda, \Lambda}}

\begin{document}
\parindent=0pt

\title[Fully nonlinear equations in radial domains]{ Existence results for fully nonlinear  equations in radial domains}

\author[ G. Galise, F. Leoni \& F. Pacella]
{Giulio Galise, Fabiana Leoni \& Filomena Pacella}

\address{Dipartimento di Matematica\newline
\indent Sapienza Universit\`a  di Roma \newline
 \indent   P.le Aldo  Moro 2, I--00185 Roma, Italy.}
 \email{galise@mat.uniroma1.it}
\email{leoni@mat.uniroma1.it}
\email{pacella@mat.uniroma1.it}

\begin{abstract}
We consider the fully nonlinear problem 
\begin{equation*}
\begin{cases}
-F(x,D^2u)=|u|^{p-1}u & \text{in $\Omega$}\\
u=0 & \text{on $\partial\Omega$}
\end{cases}
\end{equation*}
where $F$ is uniformly elliptic, $p>1$ and $\Omega$ is either an annulus or a ball in $\Rn$, $n\geq2$. \\
We prove the following results:
\begin{itemize}
	\item[i)] existence of a positive/negative radial solution for every exponent $p>1$, if $\Omega$ is an annulus;
	\item[ii)] existence of infinitely many sign changing radial solutions for every $p>1$, characterized by the number of nodal regions, if $\Omega$ is an annulus;
	\item[iii)] existence of infinitely many sign changing radial solutions characterized by the number of nodal regions, if $F$ is one of the Pucci's operator, $\Omega$ is a ball and $p$ is subcritical.
\end{itemize}
 \end{abstract}
\maketitle

\section{Introduction}\label{intro}

In this paper we study the existence of radial solutions, both positive/negative or sign changing, of the following fully nonlinear elliptic problem:
\begin{equation}\label{eq1_intro}
\begin{cases}
-F(x,D^2u)=|u|^{p-1}u & \text{in $\Omega$}\\
u=0 & \text{on $\partial\Omega$}
\end{cases}
\end{equation}
where $\Omega$ is either a ball or an annulus in $\Rn$, $n\geq2$, $p>1$ and $F$ is a real continuous function in $\Rn\times\mathcal{S}_n$, with $\mathcal{S}_n$ denoting the set of symmetric $n\times n$ matrices. We will always assume that $F$ is \emph{uniformly elliptic}, see condition \eqref{ue}, and $F(x,0)\equiv0$ in $\Rn$, so that the uniform ellipticity condition implies
\begin{equation}\label{eq2_intro}
\Mm (M)\leq F(x,M)\leq \Mp (M)\quad \forall\, x\in \R^n\, ,\ M \in {\mathcal S}_n
\end{equation}
where $\Mm$ and $\Mp$ are the the Pucci's extremal operators whose definition will be recalled in Section \ref{Settings}. Since we look for radial solution of \eqref{eq1_intro} we will also assume that $F$ is \emph{radially symmetric}, that is, if $u$ is a $C^2$ radial function then $F(x,D^2u(x))$ is a radial function too (see \eqref{rad}). Before stating our theorems let us recall some previous results.\\
As far as positive (or negative) solutions are concerned an existence result in the case when $F$ is one of the Pucci's operator, i.e. either $\Mp$ or $\Mm$, is provided in \cite{QS} for any bounded smooth domain $\Omega$. It also applies to problems with a nonlinear term $f(u)$ more general than the power $u^p$, but requires a \lq\lq subcritical\rq\rq assumption which, in the case of a power nonlinearity and positive solutions, is:
\begin{equation}\label{subcritical}
\begin{cases}
p\leq p^+:=\displaystyle \frac{\tilde n_{_+}}{\tilde n_{_+}-2} & \text{if $F=\Mp$} \medskip\\
p\leq p^-:=\displaystyle\frac{\tilde n_{_-}}{\tilde n_{_-}-2} & \text{if $F=\Mm$}
\end{cases}
\end{equation}
where the dimensional parameters $\tilde n_\pm$ are defined by
\begin{equation}\label{critical}
\tilde n_{_+}=\frac{\lambda}{\Lambda}(n-1)+1\quad\;\text{and}\quad\;\tilde n_{_-}=\frac{\Lambda}{\lambda}(n-1)+1.
\end{equation}
Note that   $\tilde n_{_+}>1$ and $\tilde n_{_-}\geq2$, because $n\geq2$ and $\Lambda\geq\lambda$. When $\tilde n_{_+}\leq2$ or $\tilde n_{_-}=2$ conditions \eqref{subcritical} reduce to $p<+\infty$.\\ 
The reason for requiring this bound on the exponent $p$ is clearly understood by the method used in \cite{QS}, which is based on a fixed point argument and relies on a-priori estimates. These estimates, in turn, are related, by a blow-up procedure, to Liouville type nonexistence results in $\Rn$ or in the half space which hold for subcritical exponents (see \cite{CL, FQ, L,QS}). It is a natural and interesting question to see whether \eqref{eq1_intro}  admits solutions beyond the critical exponents in \eqref{subcritical}, at least for some domain $\Omega$.\\ At this point a comparison with the classical semilinear case, i.e. when $F$ is the Laplacian, is in order.\\ It is well known that the semilinear problem
\begin{equation}\label{eq3_intro}
\begin{cases}
-\Delta u=|u|^{p-1}u & \text{in $\Omega$}\\
u=0 & \text{on $\partial\Omega$}
\end{cases}
\end{equation}
does not admit any solution, neither positive/negative nor sign changing, when $n\geq3$, $p\geq\frac{n+2}{n-2}$ and $\Omega$ is star-shaped. On the contrary, there exist several type of bounded sets with nontrivial topology for which solutions of \eqref{eq3_intro} exist if $p$ is the critical exponent $\frac{n+2}{n-2}$, $n\geq3$, or even if $p$ is supercritical. Typical examples of such sets are domains with holes, in particular annuli.\\
Since problem \eqref{eq3_intro} is variational the existence of solutions for it is equivalent to the existence of critical point of the associated functional in the Sobolev space $H^1_0(\Omega)$. Therefore the bound on the exponent for the existence is related to the lack of compactness for the Sobolev embeddings. A finer analysis of this phenomenon shows that some compactness is restored in domains with holes at least at certain  energy levels allowing so to prove existence results for \eqref{eq3_intro} for critical or supercritical exponents. In particular, if $\Omega$ is an annulus $A$ the compact embedding of $H^1_{0,{\rm rad}}(A)=\left\{ u \in H^1_0(A)\,:\,\text{$u$ is radial}\right\}$ into $L^p(A)$, for every $p>1$, allows to prove easily that radial solutions of \eqref{eq3_intro} exist for every exponent $p>1$.\\
Coming back to the fully nonlinear problem \eqref{eq1_intro} one could ask whether similar existence results hold for supercritical exponents, though is not clear why domains with a hole should play a particular role in a nonvariational setting. \\
Motivated by this we analyze the case when $\Omega$ is an annulus and we prove the existence of solutions for every exponent $p>1$. More precisely our first result is:
\begin{theorem}\label{teo1}
Let $A_{a,b}$ be the annulus $A_{a,b} =\{ x\in \R^n \, : a<|x|<b\}$ with $a>0$. Under the assumption  \refe{ue}-\refe{rad} on the operator $F$, problem \eqref{eq1_intro} with $\Omega= A_{a,b}$ has a positive and a negative radial solution for any $p>1$.
\end{theorem}
As far as we know this is the first existence result for \eqref{eq1_intro} with a supercritical exponent. Let us observe that the result of Theorem \ref{teo1} is new even in the case $p$ satisfies \eqref{subcritical}, i.e. is subcritical. Indeed by \cite{QS} we get the existence of a positive/negative solution only when $F$ is a Pucci's operator and, even in this case, we cannot say that the solution is radial since the Gidas-Ni-Nirenberg type symmetry result of \cite{DLS} does not hold for annular domains.\\ The proof of Theorem \ref{teo1} relies on a careful study of the associated ODE problem. We point out that this analysis would be easier if $p$ is subcritical, but requires more accuracy for other $p$'s (see Remark \ref{p*}).\\
Note that the analysis of the associated ODE problem for proving existence of radial solutions has been performed in many papers in the semilinear case, but all previous methods strongly rely on the divergence form of the operator and do not straightforwardly extend to the fully nonlinear framework.\\
As a consequence of the ODE's analysis we also get existence results for some mixed boundary value problems (Theorems \ref{mix1} and \ref{mix2}).\\
In the second part of the paper we turn to the study of sign changing solutions of \eqref{eq1_intro}. For fully nonlinear equations, for example when $F$ is a Pucci's operator, the study of sign changing solutions cannot be done in the same way as for the one sign solutions. As far as we know the only available results in this direction have been obtained in \cite{AQ} by local bifurcation. \\
This can be observed also when considering nodal eigenfunctions whose existence is only known in the radial case (\cite{EFQ,II}).\\
Again looking at the semilinear problem \eqref{eq3_intro} for which infinitely many sign changing radial solutions exist in the annulus for every $p>1$ and in the ball for any $p<\frac{n+2}{n-2}$, $n\geq3$, one could ask whether similar results hold also for the fully nonlinear problem \eqref{eq1_intro}.\\
By  \lq\lq gluing\rq\rq\, together the positive and negative solutions obtained in Theorem \ref{teo1} we prove the following result in the annulus.
\begin{theorem}\label{annulus}
Let $A_{a,b}$ be an annulus as in Theorem \ref{teo1} and assume that $F$ satisfies \eqref{ue}-\eqref{rad}. Then, for any $k\in\mathbb N$,  there exist radial solutions  $u^\pm_k$ of \eqref{eq1_intro} with $\Omega=A_{a,b}$ and finite sequences $(r^\pm_{k,j})_{j=0}^k\subset[a,b]$ such that:
\begin{equation}\label{annulus_i}
a=r^\pm_{k,0}<r^\pm_{k,1}<\ldots<r^\pm_{k,k}=b;
\tag{i}
\end{equation}
\begin{equation}\label{annulus_ii+}
(-1)^{j-1}u_k^+>0\qquad\text{in}\quad A_{r^+_{k,j-1},r^+_{k,j}}\;\;\text{for}\;\; j=1,\ldots,k;
\tag{ii$_+$}
\end{equation}
\begin{equation}\label{annulus_ii-}
(-1)^{j}u_k^->0\qquad\text{in}\quad A_{r^-_{k,j-1},r^-_{k,j}}\;\;\text{for}\;\; j=1,\ldots,k.
\tag{ii$_-$}
\end{equation}
\end{theorem}
\bigskip
Finally we assume that the domain is the ball $B_R=\left\{x\in\Rn\,:\,|x|<R\right\}$ in which case the existence and uniqueness of a positive/negative solution of \eqref{eq1_intro} when $F$ is a Pucci's operator and $p$ satisfies suitable subcritical growth requirements has been proved in \cite[Theorem 5.1]{FQ}. Using this result, Theorem \ref{teo1} and again a \lq\lq gluing\rq\rq\, procedure together with a rescaling argument we can prove the existence of infinitely many sign changing solutions in the ball.
\begin{theorem}\label{sign-changing ball}
Let $\Omega$ be the ball $B_R$, $F={\mathcal{M}^\pm_{\lambda, \Lambda}}$ and $p\leq p^-$ defined in \eqref{subcritical}. Then, for any 
$k\in\mathbb N$, there exist radial solutions  $u^\pm_k$ of \eqref{eq1_intro} and finite sequences $(r^\pm_{k,j})_{j=0}^k\subset[0,R]$ such that:
\begin{equation}\label{i}
0=r^\pm_{k,0}<r^\pm_{k,1}<\ldots<r^\pm_{k,k}=R;
\tag{i}
\end{equation}
\begin{equation}\label{ii+}
u_k^+(0)>0\quad\text{and}\quad(-1)^{j-1}u_k^+>0\quad\text{in}\quad A_{r^+_{k,j-1},r^+_{k,j}}\;\;\text{for}\;\; j=1,\ldots,k;
\tag{ii$_+$}
\end{equation}
\begin{equation}\label{ii-}
u_k^-(0)<0\quad\text{and}\quad(-1)^{j}u_k^->0\quad\text{in}\quad A_{r^-_{k,j-1},r^-_{k,j}}\;\;\text{for}\;\; j=1,\ldots,k.
\tag{ii$_-$}
\end{equation}
\end{theorem}
\bigskip
Let us point out that in the previous theorem the more restrictive upper bound $p\leq p^-$ is required even for the operator $\Mp$, whose negative solutions correspond to positive solutions of $\Mm$.\\
Note that the results of Theorems \ref{teo1}, \ref{annulus} and \ref{sign-changing ball} hold also in dimension $n=1$ in any interval in $\R$, just applying simplified versions of our proofs.

We conclude by pointing out that, as for the semilinear case, we would expect that radial positive/negative solution in the annulus should be unique. The proof of this fact is not so obvious and may require a phase plane analysis of the corresponding ODE, as in \cite{FQ}. We plan to do it in a future work.

The paper is organized as follows. In Section \ref{Settings} we describe the settings and prove some preliminary results on radial solutions. In Section \ref{exi} we prove Theorem \ref{teo1}, while Section \ref{Sign changing} is devoted to the study of sign changing radial solutions in the annulus and to the proof of Theorem \ref{annulus}. Finally in Section \ref{Ball} we consider the case of the ball and prove Theorem \ref{sign-changing ball}.

\section{Settings and preliminary results on radial solutions}\label{Settings}

We consider a continuous function
$F:\R^n\times \mathcal{S}_n\to \R$, with $\mathcal{S}_n$ denoting the set of symmetric $n\times n$ matrices equipped with the usual partial ordering
$$
M\geq N \Longleftrightarrow M-N\geq 0 \Longleftrightarrow (M-N)\xi \cdot \xi \geq 0\quad \forall\, \xi\in \R^n\, .
$$
We will always assume that $F$ is  \emph{uniformly elliptic}, that is
\begin{equation}\label{ue}
\lambda \, {\rm tr}(P)\leq F(x, M+P)-F(x, M)\leq \Lambda\, {\rm tr}(P)\, ,\quad \forall\, x\in \R^n\, ,\ M, P\in {\mathcal S}_n,\ P\geq 0\, ,
\end{equation}
for positive constants $0<\lambda\leq \Lambda$.  Moreover we require that 
\begin{equation}\label{eq F(x,0)=0}
F(x, O)\equiv 0 \quad \hbox{in } \R^n\, ,
\end{equation}
so that condition \refe{ue} implies \eqref{eq2_intro}, where $\Mm$ and $\Mp$ are the Pucci's extremal operators defined respectively as
$$
\begin{array}{c}
\displaystyle \Mm (M)= \inf_{A\in \mathcal{A}_{\lambda ,\Lambda}} {\rm tr} (AM)= \lambda \sum_{\mu_i>0} \mu_i + \Lambda \sum_{\mu_i<0} \mu_i\\[1ex]
\displaystyle
\Mp (M)= \sup_{A\in \mathcal{A}_{\lambda ,\Lambda}} {\rm tr} (AM)= \Lambda \sum_{\mu_i>0} \mu_i + \lambda \sum_{\mu_i<0} \mu_i\,.
\end{array}
$$
Here $\mathcal{A}_{\lambda ,\Lambda}=\{ A\in \mathcal{S}_n\, :\, \lambda\, I_n\leq A\leq \Lambda \, I_n\}$, $I_n$ being the identity matrix in $\mathcal{S}_n$, and $\mu_1, \ldots ,\mu_n$ being the eigenvalues of the matrix $M\in \mathcal{S}_n$. For more details about the Pucci's operators see the monograph \cite{CC}.

We further assume that the operator $F$ is \emph{radially symmetric}, that is, if $u$ is a $C^2$ radial function, then $F(x,D^2u (x))$ is a radial function.
Let us recall that for a smooth radially symmetric function $u(x)=u(|x|)$, the Hessian matrix is given by
$$
D^2u(x) = \frac{u'(|x|)}{|x|} I_n + \left( u''(|x|)- \frac{u'(|x|)}{|x|}\right) \frac{x}{|x|}\otimes \frac{x}{|x|}\, .
$$
Therefore, we require that 
for any fixed $m, d \in \R$, the function
\begin{equation}\label{rad}
x\in \R^n\setminus \{ 0\}\mapsto F\left( x, \frac{d}{|x|} I_n + \left( m-\frac{d}{|x|}\right) \frac{x}{|x|}\otimes \frac{x}{|x|}\right) \ \hbox{ is radially symmetric,}
\end{equation}
where $x\otimes x$ is the matrix whose $(i,j)-$entry is $x_i x_j$, for $i, j=1,\ldots ,n$. 

Here and in the sequel we will set, for every $(m,l)\in \R^2$ and for all $x$ such that $r=|x|>0$, 
\begin{equation}\label{effe}
\mathcal{F}(r, m,l)= F\left(x, l I_n +(m-l) \frac{x}{r}\otimes \frac{x}{r}\right)\, .
\end{equation}

Looking for positive radial solutions of problem \eqref{eq1_intro} is equivalent to prove existence of solutions for the following boundary value problems for the ODE:
\begin{equation}\label{P2}
\left\{
\begin{array}{c}
\dis -\mathcal{F} \left( r, u''(r), \frac{u'(r)}{r}\right) =u^p(r)\, , \quad \hbox{for } a<r<b\, ,\\[2ex]
u(r)>0\  \hbox{ for } a<r<b\, , \ u(a)=u(b)=0
\end{array}
\right.
\end{equation}
where $\mathcal{F} :\R^3\to \R$ is the continuous function defined by \refe{effe}.

As a consequence of the uniform ellipticity assumption \refe{ue},  the ODE in \refe{P2} can be put in normal form,  that is, there exists a continuous function $\mathcal{G}: \R^3\to \R$ such that problem \refe{P2} is equivalent to 
$$
\left\{
\begin{array}{c}
\dis u''(r)= \mathcal{G} \left( r,  \frac{u'(r)}{r}, -u^p(r)\right) \, , \quad \hbox{for } a<r<b\, ,\\[2ex]
u(r)>0\  \hbox{ for } a<r<b\, , \ u(a)=u(b)=0
\end{array}
\right.
$$
with $\mathcal{G}(r, \cdot ,\cdot)$ uniformly Lipschitz continuous in $\R^2$, 
see  \cite[Lemma 5.1]{EFQ}.

It is then natural to consider positive solutions of the initial value problem
\begin{equation}\label{P3}
\left\{
\begin{array}{c}
\dis u''(r)= \mathcal{G} \left( r,  \frac{u'(r)}{r}, -u^p(r)\right) \, , \quad \hbox{for } r>a\, ,\\[2ex]
u(a)=0\, ,\ u'(a)=\alpha
\end{array}
\right.
\end{equation}
with $\alpha >0$ to be used as a shooting parameter. By Lipschitz continuity of $\mathcal{G}(r, \cdot ,\cdot)$, problem \refe{P3} has a unique positive solution $u(r)=u(r,\alpha)$ defined on a maximal interval $[a, \rho (\alpha))$, with $\rho(\alpha)\leq +\infty$ and $u(\cdot, \alpha)$ twice continuously differentiable in $[a, \rho(\alpha))$.

Let us observe that the ODE in \refe{P3} is equivalent to the ODE in \refe{P2}, and then to
$$
-F\left( r e_1,  \frac{u'(r)}{r} I_n +\left( u''(r)-\frac{u'(r)}{r}\right) e_1\otimes e_1\right)=u^p(r)\, ,
$$
which implies, by the uniform ellipticity condition \refe{eq2_intro},
\begin{equation}\label{ineq}
\left\{
\begin{array}{c}
\dis -\Mp \left( \frac{u'(r)}{r} I_n +\left( u''(r)-\frac{u'(r)}{r}\right) e_1\otimes e_1\right) \leq u^p(r) \\ [2ex]
\dis  - \Mm \left( \frac{u'(r)}{r} I_n +\left( u''(r)-\frac{u'(r)}{r}\right) e_1\otimes e_1\right)\geq u^p(r)\, .
\end{array}
\right.
\end{equation}
Since the eigenvalues of the matrix $\frac{u'(r)}{r} I_n +\left( u''(r)-\frac{u'(r)}{r}\right) e_1\otimes e_1$ are $u''(r)$, which is simple, and $\frac{u'(r)}{r}$, which has multiplicity $n-1$, we can distinguish three different cases:
\begin{itemize}
\item[$\mathbf{C_1}$:]   $u'(r)\geq 0$  and $u''(r)\leq 0$, so that $u$ satisfies
$$
\left\{
\begin{array}{c}
\dis -\lambda\,  u''(r) -\Lambda (n-1) \frac{u'(r)}{r} \leq u^p(r)\\[1ex]
\dis -\Lambda\,  u'' (r) -\lambda (n-1) \frac{u'(r)}{r} \geq u^p(r)
\end{array}
\right.
$$
\item[$\mathbf{C_2}$:]  $u'(r)\leq 0$  and $u''(r)\leq 0$, so that $u$ satisfies
$$
\left\{
\begin{array}{c}
\dis -\lambda \left( u''(r) + (n-1) \frac{u'(r)}{r}\right) \leq u^p(r)\\[2ex]
\dis -\Lambda \left( u'' (r) + (n-1) \frac{u'(r)}{r}\right) \geq u^p(r)
\end{array}
\right.
$$
\item[$\mathbf{C_3}$:] $u'(r)\leq 0$  and $u''(r)\geq 0$, so that $u$ satisfies
$$
\left\{
\begin{array}{c}
\dis -\Lambda\,  u''(r) -\lambda (n-1) \frac{u'(r)}{r} \leq u^p(r)\\[1ex]
\dis -\lambda\,  u'' (r) -\Lambda (n-1) \frac{u'(r)}{r} \geq u^p(r)
\end{array}
\right.
$$
\end{itemize}
Clearly, the above inequalities will be used repeatedly in the sequel, since they express all the information contained in the ODE of \refe{P3}. Let us  emphasize that, by the second inequality in \refe{ineq},  $u$ cannot be at the same time convex and increasing in any interval, since $u''(r)<0$ as long as $u'(r)\geq 0$. In particular, any critical point of $u$ is a local strict maximum point for $u$.

A first easy property of the function $u(\cdot, \alpha)$ is given in the following result. Hereafter, the symbols $u',\  u''$ always mean differentiation with respect to $r$.
\begin{lemma}\label{tau}
For any $\alpha>0$,  there exists $\tau (\alpha)\in (a,\rho(\alpha))$ such that 
$$
u'(r,\alpha)>0 \ \hbox{for } r< \tau(\alpha)\, ,\ u'(\tau (\alpha),\alpha)=0\, ,\   u'(r,\alpha)<0 \ \hbox{for } r>\tau(\alpha)\, .
$$
\end{lemma}
\proof 
Let us fix $\alpha >0$ and set $u(r)=u(r, \alpha)$.
Since $u'(a)=\alpha>0$,  the second inequality in \refe{ineq} implies $u''(a)<0$. Therefore, there exists a right neighborhood  of $a$ where $u$ satisfies $\mathbf{C_1}$, and $u$ continues to satisfy $\mathbf{C_1}$ as long as $u$ increases. If $u'(r)>0$ for all $r\in [a, \rho(\alpha))$, then  $\rho(\alpha)=+\infty$ and $u$ would be concave and increasing in $[a, +\infty)$. Thus, there would exist $\lim_{r\to \infty} u(r) >0$,   $\lim_{r\to \infty} u'(r)< \alpha$ and, by property $\mathbf{C_1}$, 
$$\dis \limsup_{r\to \infty} u''(r) \leq - \frac{1}{\Lambda} \lim_{r\to \infty} u^p(r)<0\, ,
$$
a contradiction to $\lim_{r\to \infty} u'(r) < +\infty$. 

Hence, there exists $\tau(\alpha )\in (a, \rho(\alpha))$ defined as the first zero of $u'$.  As already observed, any other critical point of $u$ would be a local strict maximum point for $u$, from which it follows that $u'(r)<0$ for all $r\in (\tau (\alpha), \rho(\alpha))$.

\hfill$\Box$

In the existence and uniqueness results of the subsequent sections, a crucial role will be played by the monotonicity properties   stated in the following result. Here and in the sequel we will consider the dimension-like parameter $\tilde n_{_-}$ defined in \eqref{critical} satisfying $\tilde n_{_-}>2$, since $\tilde n_{_-}=2$ corresponds to $\Lambda=\lambda$ and $n=2$, i.e. the already known case of semilinear equations in planar domains.

\begin{proposition}\label{energies}
For every $\alpha>0$, the energy function
$$
E^\alpha_1( r) = r^{2(\tilde n_{_-}-1)} \left( \frac{u'(r,\alpha)^2}{2} +\frac{u^{p+1}(r,\alpha)}{\Lambda\, (p+1)}\right) 
$$
is monotone increasing in $[\tau(\alpha), \rho(\alpha))$, and  the energy function
$$
E^\alpha_2( r) =  \frac{u'(r,\alpha)^2}{2} +\frac{u^{p+1}(r,\alpha)}{\lambda\, (p+1)} 
$$
is monotone decreasing in $[\tau(\alpha), \rho(\alpha))$. 
\end{proposition}
\proof By Lemma \ref{tau} we know that $u(\cdot ,\alpha)$ is decreasing for $r\geq \tau(\alpha)$. Hence, $u$ satisfies either inequalities in 
$\mathbf{C_2}$ or inequalities in $\mathbf{C_3}$. By differentiating $E^\alpha_1$, it then follows that
$$
\left(E^\alpha_1\right)' (r) \geq (\tilde n_{_-}-1) r^{2 (\tilde n_{_-}-1)-1} \left[ \left(\frac{\Lambda-\lambda}{\Lambda}\right) u'(r)^2 +\frac{2}{\Lambda\, (p+1)} u^{p+1}(r)\right]
$$
in any interval on the right of $\tau(\alpha)$ where $u$ is concave, and
$$
\begin{array}{rl}
\dis\left(E^\alpha_1\right)' (r) & \dis \!\!\!\geq r^{2 (\tilde n_{_-}-1)} u'(r)\left( u''(r)+\frac{u^p(r)}{\lambda}\right) +  2 (\tilde n_{_-}-1) r^{2(\tilde n_{_-}-1)-1} \left( \frac{u'(r)^2}{2} +\frac{u^{p+1}(r)}{\Lambda\, (p+1)} \right)\\[2ex]
& \dis \!\!\! \geq 2 (\tilde n_{_-}-1) r^{2 (\tilde n_{_-}-1)-1} \frac{u^{p+1}(r)}{\Lambda\, (p+1)} 
\end{array}
$$
in any interval where $u$ is convex. In both cases, for $r\geq \tau(\alpha)$, one has
$$
\left(E^\alpha_1\right)' (r)\geq 2 (\tilde n_{_-}-1) r^{2 (\tilde n_{_-}-1)-1} \frac{u^{p+1}(r)}{\Lambda\, (p+1)}
\geq 0\, .
$$ 
Analogously, by the  inequalities in $\mathbf{C_2}$ and in $\mathbf{C_3}$, it follows that $u$ satisfies for every $r\geq \tau(\alpha)$
$$
-u''\leq \frac{u^p}{\lambda}\, ,
$$
which immediately yields $\left(E^\alpha_2\right)' (r)\leq 0$ for $r\geq \tau (\alpha)$.

\hfill$\Box$

As an easy consequence of the above proposition, we have the following
\begin{corollary}\label{zero}
If $\rho(\alpha)=+\infty$, then $\dis \lim_{r\to +\infty} u(r,\alpha)=0$.
\end{corollary}
\proof  Assume $\rho(\alpha)=+\infty$. Then, by Lemma \ref{tau}, there exists $\dis \lim_{r\to +\infty} u(r,\alpha)=c\geq 0$.  Moreover, by the monotonicity of the function 
$E^\alpha_2$ established in Proposition \ref{energies}, it follows that $u'(r,\alpha)$ has a limit as $r\to +\infty$, so that $\dis \lim_{r\to +\infty} u'(r,\alpha)=0$. Furthermore, by the inequalities in $\mathbf{C_2}$ and in $\mathbf{C_3}$, one can see that, for $r\geq \tau(\alpha)$,  $u$ always satisfies
$$
-u''\geq \frac{u^p}{\Lambda} +\frac{\Lambda}{\lambda}\,  (n-1) \frac{u'}{r}\,.
$$
Hence,
$$
0\leq \limsup_{r\to +\infty} u''\leq -\frac{c^p}{\Lambda}\, ,
$$
which gives $c=0$.

\hfill$\Box$

By the above results, we see that for any $\alpha>0$ either  $\rho(\alpha)=+\infty$ and $\lim_{r\to \infty} u(r,\alpha)=0$, or $\rho(\alpha)<+\infty$ and $u(\rho(\alpha),\alpha)=0$. Moreover,  by continuous dependence on the initial data, the function $\rho(\alpha)$ will be defined and continuous in a neighborhood  of any $\alpha>0$ where $\rho(\alpha)<+\infty$.

\section{Existence of positive (negative)  solutions in annuli}\label{exi}

This section is mainly devoted to the proof of  Theorem \ref{teo1} which, by the results of the previous section, is  reduced to show that, for any given $b>a>0$,  there exists  $\alpha>0$ such that $\rho(\alpha)=b$.

As before, we denote by $u(r, \alpha)$ the unique maximal solution of the Cauchy problem \refe{P3}, defined in the interval $[a, \rho(\alpha))$, with $\rho(\alpha)\leq +\infty$ and $u(\rho (\alpha), \alpha)=0$ if $\rho(\alpha)<+\infty$. If $\alpha>0$ is kept constant or its  value is clear from the context, we simply write $u(r)$ instead of $u(r,\alpha)$.

\begin{lemma}\label{taualfa}
For $\alpha>0$, let $\tau (\alpha)$ be as in Lemma \ref{tau}. Then:
\begin{itemize}
\item[(i)] $\dis \lim_{\alpha\to +\infty} u(\tau (\alpha), \alpha)=+\infty\, ;$

\item[(ii)] $\dis \lim_{\alpha\to +\infty} \tau (\alpha)=a \, ;$

\item[(iii)] $\dis \lim_{\alpha\searrow 0} u(\tau (\alpha), \alpha)=0\, ;$

\item[(iv)] $\dis \lim_{\alpha\searrow 0} \tau (\alpha)=+\infty \, .$
\end{itemize}
\end{lemma}

\proof  Let us denote in the following $\tau=\tau(\alpha)$.

By Lemma \ref{tau}, in the interval $[a, \tau]$, $u$ is  increasing, hence concave, and $u$ satisfies $\mathbf{C_1}$.  By using the first inequality in $\mathbf{C_1}$, it is immediate to verify that, for $\tilde n_{_-}$ as in \eqref{critical}, the energy function
$$
\mathcal{E}_1 (r)= r^{2 (\tilde n_{_-}-1)} \left( \frac{u'(r)^2}{2} +\frac{u^{p+1}(r)}{\lambda\, (p+1)}\right)
$$
satisfies $\mathcal{E}_1'(r)\geq 0$ in $[a, \tau]$. Therefore $\mathcal{E}_1(\tau)\geq \mathcal{E}_1(a)$, and this yields a first lower bound on $u(\tau, \alpha)$:
\begin{equation}\label{lbtau1}
u^{p+1}(\tau, \alpha) \geq \frac{\lambda\, (p+1)}{2} \left( \frac{a}{\tau}\right)^{2 (\tilde n_{_-}-1)} \alpha^2\, .
\end{equation}
Furthermore, writing again the first inequality in $\mathbf{C_1}$ in the form
$$
\left( r^{\tilde n_{_-}-1} u'\right)' \geq - \frac{r^{\tilde n_{_-}-1}}{\lambda} u^p
$$
and integrating from $a$ to $r\in (a, \tau]$ we get
$$
r^{\tilde n_{_-}-1} u'(r)\geq a^{\tilde n_{_-}-1} \alpha - \frac{1}{\lambda}\int_a^r s^{\tilde n_{_-}-1}u^p(s)\, ds\geq a^{\tilde n_{_-}-1} \alpha -\frac{r^{\tilde n_{_-}-1}u^p(r) (r-a)}{\lambda}\, .
$$
Dividing both sides by $r^{\tilde n_{_-}-1}$ and  integrating once again from $a$ to $\tau$, we obtain a second lower bound:
\begin{equation}\label{lbtau2}
u(\tau,\alpha) \geq \frac{a^{\tilde n_{_-}-1} \alpha}{\tilde n_{_-}-2} \left( a^{2-\tilde n_{_-}}- \tau^{2-\tilde n_{_-}} \right) -\frac{u^p(\tau, \alpha) (\tau -a)^2}{2\, \lambda} \, .
\end{equation}
On the other hand, by using the second inequality in $\mathbf{C_1}$, it is easy to check that the second energy function
$$
\mathcal{E}_2(r) = \frac{u'(r)^2}{2} +\frac{u^{p+1}(r)}{\Lambda\, (p+1)}
$$
satisfies $\mathcal{E}_2' (r)\leq 0$ in $[a,\tau]$. Hence $\mathcal{E}_2 (a)\geq \mathcal{E}_2(\tau)$, that is the upper bound
\begin{equation}\label{ubtau1}
u^{p+1}(\tau, \alpha) \leq \frac{\Lambda\, (p+1)}{2}  \alpha^2\, .
\end{equation}
Moreover, again the second inequality in $\mathbf{C_1}$ and the increasing monotonicity of $u$ yield
$$
-\left( \frac{u'(r)^2}{2}\right)' \geq \frac{u^p(r) u'(r)}{\Lambda}\, .
$$
By integrating the above inequality from $r\in [a, \tau)$ to $\tau$, we obtain
$$
u'(r) \geq \sqrt{\frac{2}{\Lambda (p+1)} \left[ u^{p+1} (\tau) -u^{p+1}(r)\right]}\, .
$$
Integrating once again in $[a, \tau]$, by the change of variable $\sigma= s/u(\tau)$, we get 
$$
\dis \sqrt{\frac{2}{\Lambda (p+1)}} (\tau -a)\leq \int_0^{u(\tau)} \frac{ds}{\sqrt{u^{p+1}(\tau)-s^{p+1}}}= \frac{1}{u^{\frac{p-1}{2}}(\tau)} \int_0^1 \frac{d\sigma }{\sqrt{1-\sigma^{p+1}}}\, .
$$
Thus, setting
$$
c_p\, := \int_0^1 \frac{ds}{\sqrt{1-\sigma^{p+1}}}\, ,
$$
we obtain a second upper bound:
\begin{equation}\label{ubtau2}
\dis u(\tau, \alpha) \leq \left( \sqrt{\frac{\Lambda\, (p+1)}{2}} \frac{c_p}{\tau -a}\right)^{\frac{2}{p-1}}\, .
\end{equation}
Putting together \refe{ubtau2} and \refe{lbtau2} we further deduce
\begin{equation}\label{lbtau3}
\dis u(\tau, \alpha) \geq \frac{a^{\tilde n_{_-}-1} \alpha}{\tilde n_{_-}-2} \left( a^{2-\tilde n_{_-}}- \tau^{2-\tilde n_{_-}} \right) 
-\left( c_p \sqrt{\frac{\Lambda \, (p+1)}{2}}\right)^{\frac{2p}{p-1}} \frac{ (\tau -a)^2}{2\, \lambda (\tau-a)^{\frac{2p}{p-1}}}\, .
\end{equation}
Now, in order to prove (i), let us argue by contradiction and assume that there exist a positive constant $M$ and a diverging sequence $\alpha_k\to +\infty$ such that 
$$
u(\tau(\alpha_k), \alpha_k)\leq M\, .
$$
From \refe{lbtau1} it then follows $\tau(\alpha_k)\to +\infty$, which in turn implies, by \refe{lbtau3},
$$
u(\tau(\alpha_k), \alpha_k)\to +\infty \, ,
$$
a contradiction. This proves (i). Hence, by \refe{ubtau2}, (ii) follows. Moreover, \refe{ubtau1}  immediately gives (iii). Finally, we observe that, by concavity, one has
$$
u(\tau , \alpha)\leq \alpha (\tau -a)< \alpha \, \tau\, .
$$
From the above inequality and  \refe{lbtau1}, it then follows
$$
\tau(\alpha)^{p+1+2(\tilde n_{_-}-1)} \geq \frac{\lambda\, (p+1)}{2 \alpha^{p-1}} a^{2(\tilde n_{_-}-1)}\, ,
$$
which proves (iv).

\hfill$\Box$

In order to prove Theorem \ref{teo1}, we further need to investigate the behavior of $\rho(\alpha)$ with respect to $\alpha$. This will be done in the next proposition by using the properties of the principal (positive) eigenvalues of the Pucci's operators.

\begin{proposition}\label{rho}
\hspace{-0.5cm}
\begin{itemize}
\item[(i)] For every $M>0$  there exists a positive constant $\delta$ depending only on $M, a,  n, p, \lambda$ and $\Lambda$ such that
$$
0<\alpha \leq M\Longrightarrow a+\delta \leq \rho(\alpha)\leq +\infty\, .
$$
\item[(ii)] For $\alpha$ sufficiently large we have $\rho(\alpha )<+\infty$ and, moreover,
$$\lim_{\alpha\to +\infty} \rho(\alpha)=a\, .
$$
\item[(iii)] Finally
$$
\lim_{\alpha\to +\infty} u'(\rho(\alpha), \alpha)=-\infty\, .
$$
\end{itemize}
\end{proposition}

\proof (i) For $0<\alpha\leq M$, let us assume $\rho=\rho(\alpha)<+\infty$. Then, by the uniform ellipticity assumption \refe{eq2_intro}, the function $u(x)=u(|x|, \alpha)$ satisfies in the annulus $A_{a, \rho}$ the eigenvalue differential inequality
$$
\left\{
\begin{array}{c}
-\Mp (D^2u) \leq \left( \max_{[a,\rho]} u\right)^{p-1} u\quad \hbox{ in } A_{a,\rho}\\[2ex]
u >0 \ \hbox{ in } A_{a,\rho}\, ,\ u=0\ \hbox{ on } \partial A_{a,\rho}
\end{array}
\right.
$$
By the characterization of the sign of the principal eigenvalues in terms of the validity of the maximum principle, see \cite{BNV, BD, BEQ, QS2}, the above inequality implies that
\begin{equation}\label{lab}
\left( \max_{[a,\rho]} u\right)^{p-1} \geq \lambda_1^+ (-\Mp, A_{a,\rho})\, ,
\end{equation}
where $\lambda_1^+ (-\Mp, A_{a,\rho})$ is the principal eigenvalue of the operator $-\Mp$ in the annulus $A_{a,\rho}$ associated with positive eigenfunctions. Now, by Lemma \ref{tau},  we know that $\max_{[a,\rho]} u= u(\tau (\alpha), \alpha)$. By using the bound \refe{ubtau1} on $u(\tau(\alpha), \alpha)$ and the scaling properties with respect to the domain of $\lambda_1^+$, we then obtain
$$
\lambda_1^+(-\Mp, A_{1, \rho/a})=a^2 \lambda_1^+(-\Mp, A_{a, \rho})\leq a^2 \left( \frac{\Lambda (p+1)}{2} \alpha^2\right)^{\frac{p-1}{p+1}}\, .
$$
Since $\lambda_1^+(-\Mp, D)\to +\infty$ as ${\rm meas}(D)\to 0$, the above inequality  shows that for $\alpha$ bounded the ratio $\rho(\alpha)/a$ keeps  bounded away from 1, that is statement (i).

(ii) By Lemma \ref{taualfa} (ii),  the statement is equivalent to prove that
$$
\lim_{\alpha\to +\infty} \frac{\rho(\alpha)}{\tau(\alpha)} =1\, .
$$
Let us  argue by contradiction and  assume that there exist  a sequence $\alpha_k\to +\infty$ and  a constant $\delta >0$, such that, setting
$\tau_k=\tau(\alpha_k)$ and $\rho_k=\rho(\alpha_k)$, with possibly $\rho_k=+\infty$, one has
$$
\rho_k > (1+\delta) \tau_k\quad \hbox{ for all } k\geq 1\, .
$$
 This means that $u_k(r)=u(r,\alpha_k)$ is strictly positive in the interval $[\tau_k, (1+\delta)\tau_k]$.
 Then, by using again the uniform ellipticity condition \refe{eq2_intro}, it follows that $u_k$, as a function of $x$,  satisfies in the annulus $A_{\tau_k,r}$ the eigenvalue differential inequality
$$
u_k>0\, ,\quad -\Mm (D^2u_k(x))\geq \left( \min_{[\tau_k, r]}u_k\right)^{p-1} u_k(x)\, ,\quad x\in A_{\tau_k,r}\, ,
$$
for every $r\in (\tau_k, (1+\delta)\tau_k]$. Denoting with $\lambda_1^+(-\Mm, A_{\tau_k,r})$  the principal eigenvalue of the operator $-\Mm$ in the domain $A_{\tau_k,r}$ associated with positive eigenfunctions,  from its  very definition  it then follows that
$$
\left( \min_{[\tau_k, r]}u_k\right)^{p-1} \leq \lambda_1^+(-\Mm, A_{\tau_k,r})\, , \quad \hbox{ for all } r\in (\tau_k, (1+\delta)\tau_k]\, .
$$
By Lemma \ref{tau}, $u_k(r)$ is monotone decreasing for $\tau_k\leq r<\rho(\alpha_k)$, so that $\min_{[\tau_k, r]}u_k =u_k(r)$. Moreover, 
 by the homogeneity and monotonicity properties of $\lambda_1^+$ with respect to the domain, one has, for every  $r\in [(1+\delta/2) \tau_k, (1+\delta) \tau_k]$,  
$$
\lambda_1^+(-\Mm, A_{\tau_k,r})= \frac{1}{\tau_k^2}  \lambda_1^+\left(-\Mm ,  A_{1, \frac{r}{\tau_k} }\right)\leq 
\frac{1}{a^2} \lambda_1^+\left(-\Mm, A_{1, 1+\frac{\delta}{2}}\right)\, .
$$
We observe that  $ \lambda_1^+\left(-\Mm, A_{1, 1+\frac{\delta}{2}}\right)$ is a positive number depending only on $n,\lambda, \Lambda$ and $\delta$, and we denote it by $C_\delta$. Summing up, we have obtained the uniform estimate
\begin{equation}\label{boundeigen}
u_k(r)^{p-1} \leq \frac{C_\delta}{a^2}\, , \quad \hbox{ for all } r\in \left[ \left(1+\frac{\delta}{2}\right) \tau_k, (1 +\delta) \tau_k\right]\, .
\end{equation}
On the other hand, by Proposition \ref{energies}, we also have
$$
E^{\alpha_k}_1 (\tau_k)\leq E^{\alpha_k}_1 (r)\, , \quad \hbox{ for all } r\in [\tau_k,\rho(\tau_k))\, ,
$$
and, therefore, for all $r\in \left[ \left(1+\frac{\delta}{2}\right) \tau_k, (1 +\delta)\tau_k \right]$, 
$$
\frac{u'_k(r)^2}{2}\geq \frac{1}{\Lambda (p+1)} \left[  \frac{u_k(\tau_k)^{p+1}}{(1+\delta)^{2(\tilde n_{_-}-1)}} -\left( \frac{C_\delta}{a^2}\right)^{\frac{p+1}{p-1}} \right] =:\, M_\delta(k)\, .
$$
We notice that, by Lemma \ref{taualfa} (i), one has 
$$\lim_{k\to +\infty} M_\delta (k)= +\infty\, .
$$ 
Recalling that $u'_k(r)\leq 0$ for $r\in [\tau_k,\rho(\alpha_k))$, we then deduce
$$
-u'_k(r) \geq \sqrt{2\, M_\delta(k)}\, , \quad \hbox{ for all } r\in \left[ \left(1+\frac{\delta}{2}\right) \tau_k , (1+\delta) \tau_k \right]\, .
$$
By integration, this implies
$$
u_k\left( \left(1+\frac{\delta}{2}\right) \tau_k\right)  \geq - \int_{\left(1+\frac{\delta}{2}\right) \tau_k}^{(1+\delta) \tau_k} u'_k(r)\, dr\geq \delta \tau_k \sqrt{\frac{M_\delta(k)}{2}}\geq \delta a \sqrt{\frac{M_\delta(k)}{2}}\, ,
$$
which is a contradiction to \refe{boundeigen} in the limit as $k\to +\infty$.

(iii) Let $\alpha$ large enough so that $\rho(\alpha)<+\infty$. By Proposition \ref{energies}, we have $E^{\alpha}_1 (\rho(\alpha))\geq E^\alpha_1 (\tau(\alpha))$. This, combined with the estimate \refe{lbtau1}, yields
$$
\rho(\alpha)^{2(\tilde n_{_-}-1)} \frac{u'(\rho(\alpha),\alpha)^2}{2} \geq \tau(\alpha)^{2(\tilde n_{_-}-1)} \frac{u(\tau(\alpha),\alpha)^{p+1}}{\Lambda (p+1)}\geq \frac{\lambda}{2 \Lambda} a^{2(\tilde n_{_-}-1)} \alpha^2\, .
$$
Since $u'(\rho(\alpha), \alpha)<0$ by Hopf Boundary Lemma, we infer
\begin{equation}\label{up}
u'(\rho(\alpha),\alpha)\leq - \sqrt{\frac{\lambda}{\Lambda }} \left( \frac{a}{\rho(\alpha)}\right)^{\tilde n_{_-}-1} \alpha\, ,
\end{equation}
and the conclusion follows from statement (ii).

\hfill$\Box$

\begin{remark}\label{unifa}
{\rm All the estimates and the convergences proved in Lemma \ref{taualfa} and Proposition \ref{rho} depend on the initial point $a$, but they are uniform with respect to it, whenever $a$ varies in a bounded interval away from zero.}

\hfill$\Box$
\end{remark}

\emph{Proof of Theorem \ref{teo1}}. Let us start with the existence of a positive solution of \eqref{eq1_intro}.  For every $\alpha>0$, let $u(r,\alpha)$ be the maximal positive solution of the Cauchy problem \refe{P3}, defined on the maximal interval $[a, \rho(\alpha))$, and satisfying $u(\rho(\alpha), \alpha)=0$ if $\rho(\alpha)<+\infty$. Since $u(a, \alpha)=0$ and $u$ is a positive radial solution of 
$$
F(x, D^2u)+u^p=0\quad \hbox{ in } A_{a, \rho(\alpha)}\, ,
$$
we need to prove only that there exists $\alpha>0$ such that $\rho(\alpha)=b$. 

Let us define the set
\begin{equation}\label{D}
D\, : = \left\{ \alpha \in (0,+\infty)\, :\, \rho (\alpha) <+\infty \right\}\, .
\end{equation}
By continuous dependence on initial data for problem \refe{P3}, the set $D$ is open and $\rho$ is a continuous function on $D$. Moreover, 
by Proposition \refe{rho} (ii), $D$ is nonempty and contains a neighborhood of $+\infty$. Let 
$(\alpha^*, +\infty)$ be the unbounded connected component of $D$, with $\alpha^*\geq 0$. If $\alpha^*=0$, then
\begin{equation}\label{rho1}
\lim_{\alpha \searrow 0} \rho(\alpha)= +\infty\, ,
\end{equation}
by Lemma \ref{taualfa} (iv), since $\rho(\alpha)>\tau(\alpha)$. If $\alpha^*>0$, then
\begin{equation}\label{rho2}
\lim_{\alpha \searrow \alpha^*} \rho(\alpha)= +\infty
\end{equation}
as well, by continuous dependence on initial data. By using again Proposition \ref{rho} (ii), we deduce that the function $\rho$ maps the interval $(\alpha^*,+\infty)$ onto the interval $(a, +\infty)$. So the existence of a positive solution of \eqref{eq1_intro} in any annulus $A_{a,b}$ is achieved.

To get a negative solution we just observe that if  $F(x,M)$ is an elliptic operator satisfying \refe{ue}-\refe{rad}, then the operator
$$
G(x, M)\, : = - F(x,-M)
$$
still satisfies \refe{ue}-\refe{rad}. So, by what we have just proved, we obtain that, for any given annulus $A_{a,b}$ and any exponent $p>1$, there exists a positive radial solution $v$ of the Dirichlet boundary value problem \eqref{eq1_intro} with $F$ replaced by $G$.
Thus, the function $u=-v$ is a negative radial solution of 
$$
\left\{
\begin{array}{c}
F(x, D^2 u)+|u|^{p-1}u =0\quad \hbox{in }A_{a,b}\\[1ex]
u<0 \  \hbox{in }A_{a,b}\, ,\ u=0 \ \hbox{on } \partial A_{a,b}\,.
\end{array}\right.
$$
The proof of Theorem \ref{teo1} is complete.

\hfill$\Box$

\begin{remark}\label{p*}
{\rm The set $D$ defined by \refe{D} clearly depends on the operator $F$, on the exponent $p$ and on the inner radius $a>0$. As a general fact, we observe that, for every uniformly elliptic operator $F$ and for any $a>0$, one has
$$
p\leq \frac{\tilde n_{_-}}{\tilde n_{_-}-2}\  \Longrightarrow \ D=(0,+\infty)\, .
$$
Indeed, for any $\alpha>0$ the solution $u=u(r, \alpha)$ of the Cauchy problem \refe{P3} is a positive supersolution of
$$
-\Mm (D^2u) \geq u^p 
$$
in the annulus $A_{a, \rho(\alpha)}$. If $\rho(\alpha)=+\infty$, then $u$ would be a positive supersolution in the exterior domain $\R^n\setminus B_a$,  but, by the results of \cite{AS}, no positive supersolution exists in exterior domains for $p$ less than or equal to the critical exponent $ \tilde n_{_-}/ \tilde n_{_-}-2$.
The optimal threshold on the exponent $p$ ensuring that $D=(0,+\infty)$ depends on the operator $F$. For instance, for the Laplace operator, i.e. for $\lambda=\Lambda$, one has $D=(0,+\infty)$ if and only if $p\leq (n+2)/(n-2)$, as it can be deduced by a phase plane analysis of the solutions of the corresponding ODE problem. A detailed study of the set $D$ in connection with the question of the uniqueness of the positive radial solution is postponed to a future work.
}

\hfill$\Box$
\end{remark}

The following remark is important for the construction of sign changing solutions of \eqref{eq1_intro} as it will be done in the next section. 

\begin{remark}\label{negative}
{\rm By applying the ODE analysis to the operator $G(x,M)=-F(x,-M)$, as indicated in the proof of Theorem \ref{teo1}, we have that for every $\beta<0$, the Cauchy problem
$$
\left\{
\begin{array}{c}
\dis u''(r)= \mathcal{G} \left( r,  \frac{u'(r)}{r}, -|u|^{p-1}u(r)\right) \, , \quad \hbox{for } r>a\, ,\\[2ex]
u(a)=0\, ,\ u'(a)=\beta
\end{array}
\right.
$$
has a unique negative maximal solution $u(r,\beta)$ defined on the interval $[a,\rho(\beta))$, satisfying $\lim_{r\to +\infty}u(r,\beta)=0$ if $\rho(\beta)=+\infty$, and $u(\rho(\beta),\beta)=0$ if $\rho(\beta)<+\infty$. Moreover, for every $\beta<0$ there exists a unique $\tau(\beta)\in (a, \rho(\beta))$ such that $u'(r,\beta)<0$ for  $r\in [a, \tau(\beta))$ and $u'(r, \beta)>0$ for $r\in (\tau(\beta), \rho(\beta))$, and one has
\begin{equation}\label{neg1}
\lim_{\beta\to -\infty} u(\tau(\beta),\beta)=-\infty\, ,\quad \lim_{\beta\to -\infty} \tau(\beta)=a\, .
\end{equation}
Furthermore, for $\beta$ sufficiently small, one has $\rho(\beta)< +\infty$ and
\begin{equation}\label{neg2}
\lim_{\beta\to -\infty} \rho(\beta)=a\, ,\quad \lim_{\beta\to -\infty} u'(\rho(\beta),\beta)=+\infty\, ,
\end{equation}
and all the above limits are uniform with respect to $a$ ranging in a compact subset of $(0,+\infty)$.}

\hfill$\Box$
\end{remark}

A direct consequence of Lemma \ref{taualfa} is the existence in an annulus $A_{a,b}$ of positive or negative solutions of the mixed boundary value problem
\begin{equation}\label{bvdn}
\begin{cases}
F(x, D^2 u)+|u|^{p-1}u=0 &\text{in $A_{a,b}$}\\
u=0 &\text{for $|x|=a$}\\
 \frac{\partial u}{\partial n} =0 & \text {for $|x|=b$.}  
\end{cases}
\end{equation}

\begin{theorem}\label{mix1}
Let $A_{a,b}$ be an annulus and $p>1$. Under the assumptions \refe{ue}-\refe{rad}, problem \refe{bvdn} has a  positive and a negative radial solution.
\end{theorem}

\proof For $a>0$ and  $p>1$  fixed,  let us consider the maximal positive solution $u(r, \alpha)$  of problem \refe{P3} defined in $[a, \rho(\alpha))$, with $\alpha>0$. Let $\tau(\alpha)$ be the zero of $u'(r,\alpha)$ given by Lemma \ref{tau}. By continuous dependence on initial data, we know that $\tau$ is a continuous function on $(0,+\infty)$. Moreover,  statements (ii) and (iv) of Lemma \ref{taualfa} imply that
$$
\tau ((0,+\infty))=(a, +\infty)\, .
$$
Therefore, for any $b>a$ there exists $\alpha \in (0,+\infty)$ such that $\tau(\alpha)=b$, and this means exactly that $u(r,\alpha)$ is a  positive radial solution of \refe{bvdn}.\\ To get a negative solution one argues as in the proof of Theorem \ref{teo1}. 

\hfill$\Box$

Finally, by means of few further considerations, we can provide an existence result also for the mixed boundary value problem symmetric with respect to problem \refe{bvdn}, namely 
\begin{equation}\label{bvnd}
\begin{cases}
F(x, D^2 u)+|u|^{p-1}u=0 &\text{in $A_{a,b}$}\\
\frac{\partial u}{\partial n} =0 & \text {for $|x|=a$}\\  
u=0 &\text{for $|x|=b$.}
\end{cases}
\end{equation}

\begin{theorem}\label{mix2}
Under the assumptions \refe{ue}-\refe{rad}, for any given $A_{a,b}$ and $p>1$, problem \refe{bvnd} has a  positive and a negative radial solution.
\end{theorem}

\proof 
We only consider positive solutions, since the negative ones can be obtained as in the previous theorem. For $a>0$ fixed and any $\gamma >0$, let us consider the initial value problem
$$
\begin{cases}
v''(r)= \mathcal{G} \left( r,  \frac{v'(r)}{r}, -v^p(r)\right) & \text{for $r>a$}\\
v(r)>0 &\text{for $r>a$}\\
v(a)=\gamma\, ,\ v'(a)=0
\end{cases}
$$
and let us denote with $v(r, \gamma)$ its maximal positive solution, defined and of class $C^2$ in $[a, \sigma (\gamma))$, for some $\sigma (\gamma)\in (a, +\infty]$. Since any critical point of $v$ in $[a, \sigma (\gamma))$ is a strict local maximum point for $v$, it follows that $v'(r,\gamma)<0$ in $(a, \sigma (\gamma))$. Hence, also arguing as in the proof of Corollary \ref{zero}, one has  $v(\sigma(\gamma), \gamma)=0$ if $\sigma(\gamma) <+\infty$ and $\lim_{r\to +\infty} v(r, \gamma)=0$ if $\sigma(\gamma)=+\infty$. Moreover, $v$ satisfies  inequalities either in $\mathbf{C_2}$
or in $\mathbf{C_3}$ in $[a, \sigma(\gamma))$.
 
Let us define the set
$$
E\, : = \{ \gamma >0\, : \sigma (\gamma) <+\infty \}\, .
$$
By  the continuous dependence on initial data,  $E$ is an open subset of $(0, +\infty)$ and $\sigma$ is a continuous function on $E$. Moreover, 
by repeating words by words the proof of statement (ii) of Proposition \ref{rho}, where we used only the properties of the solution $u$ for $r>\tau$, it follows that $E$ contains a neighborhood of $+\infty$ and
\begin{equation}\label{sigma}
\lim_{\gamma\to +\infty} \sigma (\gamma) =a\, .
\end{equation}
Let $(\gamma^* ,+\infty)$ be the unbounded connected component of $E$, with $\gamma^*\geq 0$, and let us show that
\begin{equation}\label{sigma1}
\lim_{\gamma\searrow\gamma^*} \sigma (\gamma) =+\infty\, .
\end{equation}
The limit \refe{sigma1} is an easy consequence of the continuous dependence on initial data in the case $\gamma^*>0$, so let us consider only the case $\gamma^*=0$. We observe that, by the 
 increasing monotonicity of the energy function $E_1$ given by Proposition \ref{energies}, for every $\gamma\in E$ one has
\begin{equation}\label{gamma}
a^{2(\tilde n_{_-}-1)}\frac{\gamma^{p+1}}{\Lambda (p+1)} \leq \sigma(\gamma)^{2(\tilde n_{_-}-1)}\frac{v'(\sigma(\gamma),\gamma)^2}{2}\, .
\end{equation}
Furthermore, since $v$ is of class $C^2$ in $[a, \sigma(\gamma)]$ for $\gamma\in E$ and it satisfies inequalities either in $\mathbf{C_2}$
or in $\mathbf{C_3}$ in $[a, \sigma(\gamma)]$,
 we get that $v''(a, \gamma)<0$ and $v''(\sigma(\gamma),\gamma)>0$. Thus, there exists $\sigma_1\in (a, \sigma(\gamma))$ such that 
 $$
 v''(\sigma_1, \gamma)=0\ \hbox{ and } v''(r, \gamma)>0\ \hbox{ for } r\in (\sigma_1, \sigma(\gamma)]\, .
 $$
Hence, by using also the first inequality of $\mathbf{C_3}$ evaluated at $\sigma_1$, we obtain
$$
v'(\sigma(\gamma), \gamma) > v'(\sigma_1, \gamma) \geq - \frac{ \sigma_1 \, v(\sigma_1, \gamma)^p}{\lambda (n-1)}\, .
$$ 
By observing that $v'(\sigma(\gamma), \gamma)<0$, $v(\sigma_1, \gamma)< \gamma$ and $\sigma_1 <\sigma(\gamma)$, it then follows
$$
v'(\sigma(\gamma), \gamma)^2 < \frac{ \sigma (\gamma)^2 \, \gamma^{2p}}{\lambda^2(n-1)^2}\, ,
$$
which yields, together with \refe{gamma},
$$
\frac{2 \left( a^{(\tilde{n}-1)}\lambda (n-1)\right)^2}{\Lambda (p+1) \gamma^{p-1}} < \sigma(\gamma) ^{2 \tilde n_{_-}}\, .
$$
This proves  that $\sigma(\gamma)\to +\infty$ as $\gamma\to 0$, that is \refe{sigma1} in the case $\gamma^*=0$. Hence, by \refe{sigma} and the continuity of $\sigma$ on $E$, it follows that
$$
\sigma\left( (\gamma^*, +\infty)\right) = (a, +\infty)\, ,
$$
that is, for every $b>a$ there exists $\gamma\in E$ such that $\sigma (\gamma)=b$. This means exactly that $v=v(r, \gamma)$ is a radial solution of problem \refe{bvnd}.

\hfill$\Box$

\section{Sign changing solutions in annuli}\label{Sign changing}
This section is devoted to the proof of Theorem \ref{annulus}. For any $k\in\mathbb N$ we aim to construct classical solutions $u_k(x)=u_k(r)$ of \eqref{eq1_intro} whose nodal sets ${\mathcal N}_{u_k}=\left\{x\in A_{a,b}:\,u_k(x)=0\right\}$ consist exactly of $k-1$ concentric spheres.

With the same notations of Section \ref{Settings}, for $\alpha>0$ let us consider the initial value problem 
\begin{equation}\label{ivp}
\left\{
\begin{array}{rl}
u''(r)= \mathcal{G} \left(r,\frac{u'(r)}{r}, -|u|^{p-1}u(r)\right) & \text{for $r>a$}\\
u(a)=0,\,\;u'(a)=\alpha.&
\end{array}\right.
\end{equation}
We will show that, if $\alpha$ is suitably chosen, then \eqref{ivp} admits solutions with any prescribed number of zeros and satisfying further the Dirichlet boundary condition $u(b)$=0.\\ Using the results of the previous sections concerning the existence of positive/negative solutions of  \eqref{eq1_intro} and their qualitative properties, we first construct oscillating solutions of problem \eqref{ivp} for large $\alpha$.

\begin{proposition}\label{nodal}
For any $k\in\mathbb N$ there exists a nonnegative constant $\alpha^*_k$ such that for any $\alpha>\alpha_k^*$ there exists a finite sequence
\begin{equation}\label{i_prop_annulus}
a=r^+_{k,0}(\alpha)<\ldots<r^+_{k,k}(\alpha)
\tag{i}
\end{equation}
and a corresponding solution $u^+_k(r)=u^+_k(r,\alpha)$ of \eqref{ivp} in $[a,r^+_{k,k}(\alpha)]$, satisfying the following properties:
\begin{equation}\label{ii_prop_annulus}
\lim_{\alpha\to+\infty}r^+_{k,k}(\alpha)=a\quad\text{and}\quad \lim_{\alpha\searrow\alpha^*_k}r^+_{k,k}(\alpha)=+\infty;
\tag{ii}
\end{equation}
\begin{equation}\label{iii_prop_annulus}
u^+_k(r^+_{k,j}(\alpha))=0 \quad\text{and\quad $(-1)^{j-1}u^+_k(r)>0$\quad in\quad$(r^+_{k,j-1}(\alpha),r^+_{k,j}(\alpha))$\quad for \quad$j=1,\ldots,k$};
\tag{iii}
\end{equation}
\begin{equation}\label{iv_prop_annulus}
\lim_{\alpha\to+\infty}(-1)^{j}(u^+_k)'(r^+_{k,j}(\alpha))=+\infty\quad\text{for \quad$j=0,\ldots,k$}.
\tag{iv}
\end{equation}
Moreover for any $j=1,\ldots,k$ the mapping $\alpha\mapsto r^+_{k,j}(\alpha)$ is continuous in $(\alpha_k^*,+\infty)$.
\end{proposition}

\proof
We argue by induction on $k\in\mathbb N$. In the case $k=1$, the existence of $\alpha^*_1$, $r^+_{1,1}$ and $u^+_1$ just follows from Proposition \ref{rho}. We let $u^+_1(r)=u(r,\alpha)$ be the solution of problem \eqref{ivp} which is positive in the interval $(a,\rho(\alpha))$, with $\rho(\alpha)$ satisfying (ii) of Proposition \ref{rho}. Therefore, it suffices to define  $\alpha^*_1:=\alpha^*$ as the infimum of the unbounded connected component of the set $D_1=\left\{ \alpha \in (0,+\infty)\, :\, \rho (\alpha) <+\infty \right\}$, and 
$r^+_{1,1}:=\rho(\alpha)$. Clearly \eqref{i_prop_annulus} is satisfied. Condition \eqref{ii_prop_annulus} follows from Proposition \ref{rho} (ii) and from \eqref{rho1} and \eqref{rho2}; \eqref{iii_prop_annulus}  holds true by the very definition of $u^+_1$, and \eqref{iv_prop_annulus} is just Proposition \ref{rho} (iii). The continuity of $\alpha\mapsto r_{1,1}^+(\alpha)$ follows from the continuous dependence on initial data in \eqref{ivp}.\\
Let us now assume that all the statements are true for $k\geq1$, and let us prove them for $k+1$. For $\alpha>\alpha^*_k$, we extend the function $u^+_k$ by considering the solution of the Cauchy problem
\begin{equation*}
\left\{
\begin{array}{rl}
v_{k+1}''(r)= \mathcal{G} \left(r,\frac{v_{k+1}'(r)}{r}, -|v_{k+1}|^{p-1}v_{k+1}(r)\right) & \text{for $r> r^+_{k,k}(\alpha)$}\\
v_{k+1}(r^+_{k,k}(\alpha))=0,\,\;v_{k+1}'(r^+_{k,k}(\alpha))=(u^+_k)'(r^+_{k,k}(\alpha)) &
\end{array}\right.
\end{equation*}
and we define $\rho_{k+1}(\alpha)\leq+\infty$ as the end point of the maximal interval $[r^+_{k,k}(\alpha),\rho_{k+1}(\alpha)]$ such that 
$$
(-1)^k v_{k+1}>0\quad \text{in}\quad(r^+_{k,k}(\alpha),\rho_{k+1}(\alpha)).
$$
By the induction assumption, Proposition \ref{rho} and Remarks \ref{unifa} and \ref{negative}, it follows that $\rho_{k+1}(\alpha)<+\infty$ for $\alpha$ sufficiently large and $\lim_{\alpha\to+\infty}\rho_{k+1}(\alpha)=a$. Thus we consider the nonempty open set $D_{k+1}=\left\{ \alpha\in (\alpha^*_k ,+\infty)\, :\, \rho_{k+1}(\alpha) <+\infty \right\}$ and we define $\alpha^*_{k+1}\geq\alpha^*_k$ as the infimum of the unbounded connected component of $D_{k+1}$. For $\alpha>\alpha^*_{k+1}$, we then define  
\begin{equation*}
r^+_{k+1,j}(\alpha):=
\begin{cases}
r^+_{k,j}(\alpha) & \text{for $j=0,\ldots,k$}\\
\rho_{k+1}(\alpha) & \text{for $j=k+1$}
\end{cases}
\end{equation*}
and 
$$
u^+_{k+1}(r):=
\begin{cases}
u^+_k(r) & \text{if $r\in[a,r^+_{k+1,k}(\alpha)]$}\\
v_{k+1}(r) & \text{if $r\in(r^+_{k+1,k}(\alpha),r^+_{k+1,k+1}(\alpha)]$}.
\end{cases}
$$
Then \eqref{i_prop_annulus} is satisfied by construction. Moreover $u^+_{k+1}$ clearly is a  solution of \eqref{ivp} and statements \eqref{ii_prop_annulus}, \eqref{iii_prop_annulus} and \eqref{iv_prop_annulus} follow by using the properties of $r^+_{k,k}$ and by applying again Proposition \ref{rho}  and Remarks \ref{unifa} and \ref{negative}.

\hfill$\Box$

As a consequence of the previous proposition we easily obtain the assertion of Theorem \ref{annulus}.

\bigskip

\emph{Proof of Theorem \ref{annulus}}. It is sufficient to focus on $u_k^+$, since the case $u_k^-$ can be easily converted into the previous one as in the proof of Theorem \ref{teo1}.\\
In view of Proposition \ref{nodal}, for each $k\in\mathbb N$ and $\alpha$ large there exists a nodal solution $u^+_k(x)=u^+_k(|x|,\alpha)$ of the problem
\begin{equation*}
\begin{cases}
-F(x, D^2 u)=|u|^{p-1}u  & \text{in $A_{a,r^+_{k,k}(\alpha)}$}\\
 u=0 & \text{on $\partial A_{a,r^+_{k,k}(\alpha)}$}\,.
\end{cases}
\end{equation*}
Since the mapping $\alpha\mapsto r^+_{k,k}(\alpha)$ is continuous in a neighborhood of infinity with range $(a,+\infty)$, we can pick $\alpha=\alpha(b)$ in order to adjust the Dirichlet boundary condition $r^+_{k,k}(\alpha)=b$. 

\hfill$\Box$

\section{Radial solutions in balls}\label{Ball}

In this last section, we address the issue of the existence  sign changing radial solutions of the  problems 
\begin{equation}\label{Lane-Emden Pucci_ball}
\begin{cases}
-{\mathcal{M}^\pm_{\lambda, \Lambda}}(D^2u)=|u|^{p-1}u  & \text{in $B_R$}\\
 u=0 & \text{on $\partial B_R$},
\end{cases}
\end{equation}
where $B_R$ is the ball of radius $R>0$ centered at the origin and ${\mathcal{M}^\pm_{\lambda, \Lambda}}$ means either $\Mp$ or $\Mm$. Contrary to the case of annular domains, now the existence of nontrivial solutions strongly depends on the exponent $p>1$ of the power nonlinearity. As far as the well-posedness of \eqref{Lane-Emden Pucci_ball} is concerned, in the class of positive functions we quote the following result (see \cite[Theorem 5.1]{FQ}).

\begin{theorem}\label{uniq_ball} There exist critical exponents $p^*_\pm>1$ such that each of the problems
\begin{equation}\label{Lane-Emden Pucci_ball_positive}
\begin{cases}
-{\mathcal{M}^\pm_{\lambda, \Lambda}}(D^2u)=u^p  & \text{in $B_R$}\\
 u>0 & \text{in $B_R$}\\
 u=0& \text{on $\partial B_R$}
\end{cases}
\end{equation}
has a unique solution if, and only if, $p<p^*_\pm$ .
\end{theorem}

\begin{remark}
\rm Note that, by the symmetry results of \cite{DLS}, the solutions of \eqref{Lane-Emden Pucci_ball_positive} are radial.
The  exponents $p^*_\pm$ of the previous theorem are respectively strictly bigger than $p^\pm$, defined in \eqref{subcritical}. As further proved in \cite{FQ}, the numbers $p^*_\pm$ act as critical exponents for the existence of positive radial solutions of ${\mathcal{M}^\pm_{\lambda, \Lambda}}(D^2u)+u^p=0$ in $\mathbb R^n$. 

\hfill$\Box$
\end{remark}

In order to build nodal solutions of the problems \eqref{Lane-Emden Pucci_ball}, we will use  a gluing procedure between the solutions with constant sign in annuli, provided by Theorem \ref{teo1}, and those of Theorem \ref{uniq_ball} in the ball, taking advantage of the positive homogeneity of the Pucci's extremal operators. In particular we recall that if $u$ is a solution of $$-{\mathcal{M}^\pm_{\lambda, \Lambda}}(D^2u)=|u|^{p-1}u\quad\text{ in $B_R$,}$$ then  for any $\varrho>0$ the function $v(x)=\varrho^{\frac{2}{p-1}}u(\varrho x)$ is in turn a solution in $B_{\frac R\varrho}$. The proof of thi fact is a straightforward computation and it is still true for sub/supersolutions.\\ The argument we propose also rely on the availability of Liouville type theorems. For this reason we further restrict our framework to the subcritical assumption $p\leq p^-$.  
As far as Liouville type results are concerned, we refer to \cite{AS} for  the case of exterior domains, and to \cite{CL} and \cite{L}  for the whole space $\mathbb R^n$ and the half space $\partial \mathbb R^n_+$ respectively.

\bigskip
 
\emph{Proof of Theorem \ref{sign-changing ball}}. As in the proof of Theorem \ref{annulus} we only treat the case $u^+_k$ and without loss of generality we may consider the operator $\Mp$. We proceed by induction on $k\in\mathbb N$, observing that the initial step $k=1$ is a direct consequence of Theorem \ref{uniq_ball}. Let us now assume that $u^+_k(x)=u^+_k(r)$ satisfies \eqref{i}-\eqref{ii+} for $k\geq1$ and extend it by considering the function $v_{k+1}(x)=v_{k+1}(r)$, solution  of the initial value problem 
\begin{equation*}
\left\{
\begin{array}{rl}
v_{k+1}''(r)= \mathcal{G} \left(\frac{v_{k+1}'(r)}{r}, -|v_{k+1}|^{p-1}v_{k+1}(r)\right) & \text{for $r\geq R$}\\
v_{k+1}(R)=0,\,\;v_{k+1}'(R)=(u^+_k)'(R) &
\end{array}\right.
\end{equation*}
defined in the maximal interval $[R,\rho_{k+1}]$, with $\rho_{k+1}$ is such that
\begin{equation}\label{eq_v_k+1}
(-1)^kv_{k+1}>0\qquad\text{in}\quad A_{R,\rho_{k+1}}.
\end{equation}
The existence of $\rho_{k+1}$ follows from the positivity of the normal derivative $(-1)^k(u^+_k)'(R)$, as a consequence of Hopf Boundary Lemma. Moreover $\rho_{k+1}<+\infty$, since otherwise the function $v(x):=(-1)^kv_{k+1}(x)$ would be a positive supersolution  of 
$$
-\Mm(D^2v)\geq v^p\qquad\text{in}\quad \mathbb R^n\backslash B_{R}
$$
in the subcritical case $1<p\leq p^-$, so contradicting the nonexistence result \cite[Theorem 1.4 (i)]{AS} of nontrivial nonnegative supersolutions in exterior domains.\\
We now consider the function
$$
\tilde u^+_{k+1}(x)=\tilde u^+_{k+1}(|x|):=
\begin{cases}
u^+_k(|x|) & \text{if $x\in \overline B_{R}$}\\
v_{k+1}(|x|) & \text{if $x\in \overline A_{R,\rho_{k+1}}$}. 
\end{cases}
$$
By construction $\tilde u^+_{k+1}$ is a solution of \eqref{Lane-Emden Pucci_ball} in the ball $B_{\rho_{k+1}}$. In order to produce a nodal solution in $B_R$ we only need to rescale this function by defining
$$
u_{k+1}^+(x):=\left(\frac{\rho_{k+1}}{R}\right)^{\frac{2}{p-1}}\tilde u^+_{k+1}\left(\frac{\rho_{k+1}}{R}x\right)
$$
and the corresponding finite sequence of its zeros
$$
r_{k+1,j}:=\begin{cases}
\frac{R}{\rho_{k+1}} r^+_{k,j} & \text{for $j=1,\ldots,k$}\\
R & \text{for $j=k+1$}.
\end{cases}
$$
The function $u^+_{k+1}$ satisfies the boundary condition in $B_R$. Moreover, condition \eqref{i} works at level $k+1$ by the induction  scheme and the fact that $R<\rho_{k+1}$. Concerning \eqref{ii+} we first note that
$$
u^+_{k+1}(0)=\left(\frac{\rho_{k+1}}{R}\right)^{\frac{2}{p-1}}u^+_{k}(0)>0.
$$
Furthermore if $j=1,\ldots,k$ and $x\in A_{r^+_{k+1,j-1},{r^+_{k+1,j}}}$, then $\frac{\rho_{k+1}}{R}x\in A_{r^+_{k,j-1},r^+_{k,j}}$ and
$$
(-1)^{j-1}u^+_{k+1}(x)=\left(\frac{\rho_{k+1}}{R}\right)^{\frac{2}{p-1}}(-1)^{j-1}\tilde u^+_{k+1}\left(\frac{\rho_{k+1}}{R}x\right)=\left(\frac{\rho_{k+1}}{R}\right)^{\frac{2}{p-1}}(-1)^{j-1}u^+_{k}\left(\frac{\rho_{k+1}}{R}x\right)>0,
$$
again by the induction scheme. Finally for $j=k+1$ and $x\in A_{r^+_{k+1,k},r^+_{k+1,k+1}}=A_{\frac{R^2}{\rho_{k+1}},R}$ we have
$$
(-1)^ku^+_{k+1}(x)=\left(\frac{\rho_{k+1}}{R}\right)^{\frac{2}{p-1}}(-1)^{k}\tilde u^+_{k+1}\left(\frac{\rho_{k+1}}{R}x\right)=\left(\frac{\rho_{k+1}}{R}\right)^{\frac{2}{p-1}}(-1)^{k}v_{k+1}\left(\frac{\rho_{k+1}}{R}x\right)>0,
$$
because  $\frac{\rho_{k+1}}{R}x\in  A_{R,\rho_{k+1}}$ and \eqref{eq_v_k+1}. The proof is complete.

\hfill$\Box$

\end{document}